\documentclass[11pt]{amsart}

\usepackage[height=190mm,width=140mm]{geometry}
\usepackage{amscd,amssymb,amsfonts,amsmath,amsthm}
\usepackage{verbatim}
\usepackage{fancyheadings}
\usepackage{xspace}
\usepackage{pstricks}
\usepackage[ansinew]{inputenc}
\usepackage[all]{xy}

\usepackage{hyperref}

\makeindex

\newcommand{\R}{\mathbb{R\/}}

\newcommand{\tK}{\tilde{K}}
\newcommand{\tF}{\tilde{F}}
\newcommand{\tR}{\tilde{R}}
\newcommand{\tE}{\tilde{E}}

\newcommand{\Kb}{\overline{K}}

\newcommand{\Kt}{K(t)}

\newcommand{\C}{\mathbb{C\/}}
\newcommand{\Q}{\mathbb{Q\/}}

\newcommand{\Ct}{{\mathbb{C\/}(t)}}
\newcommand{\Rt}{{\mathbb{R\/}(t)}}
\newcommand{\Cz}{{\mathbb{C\/}(z)}}
\newcommand{\Rz}{{\mathbb{R\/}(z)}}

\def\P{\mathbb{P}^1}

\newcommand{\Gm}{ \mathbb{G\/}_m}

\newcommand{\m}{\mathfrak m\/}

\newcommand{\D}{ {\sf DRing}}

\newcommand{\PV}{ {\sf PV}}

\newcommand{\lra}{\longrightarrow}
\newcommand{\hra}{\hookrightarrow}

\newcommand{\PVR}{Picard-Vessiot ring\xspace}

\newcommand{\PVF}{Picard-Vessiot field\xspace}

\newcommand{\PVE}{Picard-Vessiot extension\xspace}

\def\O{\mathcal{O\/}}

\DeclareMathOperator{\Gal}{\underline{Gal}^\partial}

\DeclareMathOperator{\SO}{SO}
\DeclareMathOperator{\GL}{GL}
\DeclareMathOperator{\M}{Mat}

\DeclareMathOperator{\Hom}{Hom}
\DeclareMathOperator{\Aut}{Aut}
\DeclareMathOperator{\Isom}{Isom}
\DeclareMathOperator{\Heins}{H^1}

\DeclareMathOperator{\Quot}{Quot}
\DeclareMathOperator{\Spec}{Spec}
\DeclareMathOperator{\id}{id}

\newcommand{\pd}[2]{\frac{\partial #1}{\partial #2}}

\newtheorem{theo}{Theorem}[section]

\newtheorem{lem}[theo]{Lemma}
\newtheorem{cor}[theo]{Corollary}

\newtheorem{prop}[theo]{Proposition}
\newtheorem{defi}[theo]{Definition}

\newcommand{\defind}[1]{\index{#1}{\bf #1}}

\begin{document}

\title[The real inverse problem]{The inverse problem of differential Galois theory over the field
$\Rz$}

\author{Tobias Dyckerhoff}
\address{ University of Pennsylvania\\ Department of Mathematics\\ David Rittenhouse Laboratory\\
Philadelphia, PA 19104}

\email{tdyckerh@math.upenn.edu}

\begin{abstract} We describe a Picard-Vessiot theory for differential fields 
	with non algebraically closed fields of constants. As a technique for constructing and
	classifying Picard-Vessiot extensions, we develop a Galois descent theory.
	We utilize this theory to prove that every linear algebraic group $G$ over $\R$ occurs as a
	differential Galois group over $\Rz$. The main ingredient of the proof is the Riemann-Hilbert
	correspondence for regular singular differential equations over $\Cz$. 
\end{abstract}

\maketitle

\section{Introduction}

Given a differential field $F$ with field of constants $K$ and a linear algebraic group $G$ over
$K$, does there exist a field extension $E/F$ with differential Galois group $G$?\\
This problem is known as the inverse problem of differential Galois theory. For the classical case $F =
\Cz$, the first complete answer was given by C. Tretkoff and M. Tretkoff in 1979:

\begin{theo}
Given any linear algebraic group $G$ over $\C$, there exists a Picard-Vessiot extension $E$ of $\Cz$
whose differential Galois group is isomorphic to $G$.
\end{theo}
\begin{proof}\cite{tretkoff}\end{proof}

The main ingredient of the proof is the Riemann-Hilbert correspondence for regular singular
differential equations on $\P$ which is a transcendental result.\\
Since then, a positive answer to the inverse problem has been obtained over $K(z)$ for any
algebraically closed field $K$ (cf. \cite{mitschi}, \cite{hartmann}). The proof is purely algebraic
and uses the structure theory of linear algebraic groups.\\

In this paper, we generalize the classical result to the following theorem.
\begin{theo}
Given any linear algebraic group $G$ over $\R$, there exists a Picard-Vessiot extension $E$ of $\Rz$
whose differential Galois group is isomorphic to $G$.
\end{theo}

We remark that this theorem also generalizes the inverse problem of finite Galois theory over $\Rz$
which was solved by W. Krull and J. Neukirch in \cite{neukirch1}.\\

To prove the theorem, we first establish an appropriate Picard-Vessiot theory over non algebraically 
closed fields of constants of characteristic zero. 
This is done in Section \ref{basicpv}. We use ideas from \cite{levelt}
as a starting point, but we give simplifications and additional results. Given a differential
equation over a differential field $F$, there are two immediate natural questions: Does a
Picard-Vessiot extension exist and if so, is it unique? Both questions have a positive answer in the
case of an algebraically closed field of constants. We are able to prove the following existence
result in the general case:
\begin{theo}
	Let $S$ be a $\partial$-ring which is noetherian, integral and of finite Krull dimension. 
	Let $F$ be the field of fractions of $S$ whose field of constants we denote by $K$.
	Assume that $\Spec(S)$ has a $K$-valued point which is regular. 
	Then for any differential equation with coefficients in $S \subset F$ there exists a \PVE $E/F$.
\end{theo}
The uniqueness of a \PVE fails, but we obtain a classification result which will be explained in the
following paragraph.\\

In Section \ref{descent}, we develop a descent theory for Picard-Vessiot extensions. This is an
essential tool which we will use to construct and classify such extensions.\\
Let $F$ is a differential field
with field of constants $K$ and $L/K$ an algebraic Galois extension. Then our main result is an equivalence
between the category of Picard-Vessiot extensions $E/F$ and such extensions over $F\otimes_K L$
equipped with an action of the Galois group of $L/K$ (Theorem \ref{pve.descent}).
An immediate consequence of this theorem is a classification of Picard-Vessiot extensions
by Galois cohomology:
Let $E/F$ and $E'/F$ be Picard-Vessiot extensions and $K$ be the field of constants of $F$. 
For a field extension $L/K$
we call $E'$ an \defind{$L/K$-form} of $E$ if $E_L$ and $E'_L$ are isomorphic in $(\D/F_L)$.
	
\begin{cor}
	Let $E/F$ be a \PVE with $\partial$-Galois group $G$ and let $L/K$ be an algebraic 
	Galois extension with Galois group $\Gamma$. 
	Then the pointed Galois cohomology set $\Heins(\Gamma, G(L))$ classifies isomorphism classes of
	$L/K$-forms of $E/F$.
\end{cor}

The succeeding section contains a geometric proof of the differential Galois Correspondence in the case of 
a non algebraically closed field of constants. We do not refer to the descent theory from Section
\ref{descent} which would yield a reduction to the algebraically closed case, but prefer to give a
direct proof using Grothendieck's theory of descent for quasi-projective schemes \cite{sga1}.\\
We remark that, restricted to the case of finite Galois extensions, our theory is slightly more
general than the usual finite Galois theory. Namely, we extend the correspondence to 
finite \'etale $G$-torsors for a finite group scheme $G$ over $K$. These are classically Galois if
and only if $G(K) = G(\Kb)$.\\

Section \ref{sec.ex} illustrates some phenomena of the developed theory by means of examples.\\

In Section \ref{monodromy} we establish some results on the interaction of the complex conjugation 
in the complex plane and the monodromy representation of a singular differential equation. The key
proposition which provides the link between the descent theory and the Riemann-Hilbert
correspondence is proved in Section \ref{riemann}:
\begin{prop}\label{supi}
	Let $\tilde{E}/\Cz$ be a \PVE for the matrix $A \in \M_n(\Cz)$. 
	If the equation $\pd{y}{z} = Ay$ is equivalent to its conjugate
	$\pd{y}{z} = (\tau.A) y$ then $\tilde{E}/\Cz$ descends to a \PVE $E/\Rz$.
\end{prop}

Combining the developed techniques, we prove the main theorem in Section \ref{finalsection}. Roughly, 
by choosing an appropriate monodromy representation, we can assure that the corresponding
differential equation will be equivalent to its complex conjugate. This allows us to apply
Proposition \ref{supi} and obtain the final result.

\section{Basic Picard-Vessiot theory}\label{basicpv}

In this section we establish basic results of Picard-Vessiot theory over non algebraically
closed fields of constants $K$. We mainly follow \cite{levelt} but since we have some simplifications and
additional results we will give detailed proofs.\\

We introduce some conventions. It is understood that all fields under consideration are of
characteristic zero. Whenever the symbol $\partial$ occurs in connection with an English word it
should be read \emph{differential}. All derivations will be denoted by $\partial$. Derivations are
extended to tensor products via the \defind{Leibniz rule}
\[
\partial(x \otimes y) = \partial(x) \otimes y + x \otimes \partial(y)
\]
and to localizations via the \defind{quotient rule}
\[
\partial \left( \frac{x}{y} \right) = \frac{y\partial(x) - x \partial(y)}{y^2} \text{.}
\]
We work over a field $K$ of characteristic zero which is considered as a differential field
with the derivation $\partial=0$.\\

Let $\D$ be the category of $\partial$-rings. The objects are pairs $R = (R,\partial)$ where $R$ is
a commutative ring with unit and $\partial$ is a derivation on $R$. The morphisms are ring
homomorphisms $f: R \to S$ with $f\circ \partial = \partial\circ f$.\\
Given a $\partial$-ring $S$ we define the category $(\D/S)$ whose objects are morphisms $S \to R$ in $\D$ 
and whose morphisms are commutative diagrams
\[
\xymatrix@C-1.2pc{R \ar[rr] & & R'\\ & S\text{.}\ar[ul]\ar[ur]&}
\]
For a morphism $S \to S'$ in $\D$ we define the base change 
functor $-_{S'} : (\D/S) \to (\D/S')$ mapping $R$ to $R \otimes_S S'$.\\
We define an \defind{extension of $\partial$-rings} $R/S$ to be an injection $S \hra R$ in $\D$. 
An extension of $\partial$-rings $R/S$ in $(\D/K)$ is called \defind{geometric} 
if $R$ and $S$ are integral domains whose fields of fractions have $K$ as field of constants.
The field $K$ is then called the \defind{field of constants} of $R/S$.
A $\partial$-ring $S$ is called \defind{$\partial$-simple} if it has no proper, 
nontrivial $\partial$-ideals.\\ 

Let $F$ be a differential field of characteristic zero with derivation $\partial$ and \defind{field of
constants} $K = \{ f \in F\, |\: \partial(f) = 0 \}$. 
Consider the linear differential equation 
\[
\partial\left( \begin{array}{c} w_1\\ \vdots\\ w_n \end{array}\right) = A\left( \begin{array}{c} w_1\\ \vdots\\ w_n \end{array}\right)
\]
defined by the $n$ by $n$ matrix $A \in \M_n(F)$.
A \defind{Picard-Vessiot extension} $E/F$ for this equation is defined to be a differential field over $F$
with the properties
\begin{enumerate}
	\item $E/F$ is generated by a fundamental solution matrix for $A$, i.e. there exists a matrix
		$W \in \GL_n(E)$ such that $\partial(W) = AW$ and $E=F(W_{ij})$.
	\item The field of constants of $E$ is $K$.
\end{enumerate}

To a Picard-Vessiot extension $E/F$ we associate the group functor
\[
\underline{\Aut}^\partial(E/F) : ({\sf Algebras }/K) \lra ({\sf Groups}),\; 
L \mapsto \Aut^\partial(E\otimes_K L / F \otimes_K L)
\]
which will turn out to be representable by an affine scheme of finite type over $K$. By Yoneda's lemma
we conclude that the latter is a linear algebraic group over $K$.\\

It turns out to be fruitful to replace the Picard-Vessiot extension $E/F$ by a more geometric
object, the Picard-Vessiot ring. 

\begin{defi}\label{dgalois}
Consider the linear differential equation $\partial(w) = Aw$ defined by the $n$ by $n$ matrix $A \in \M_n(F)$.
A \defind{Picard-Vessiot ring} $R/F$ for this equation is defined to be a $\partial$-ring over $F$
having the properties
\begin{enumerate}
	\item $R/F$ is generated by a fundamental solution matrix for $A$, i.e. there exists a matrix
		$W \in \GL_n(R)$ such that $\partial(W) = AW$ and $R=F[W_{ij}, \det(W)^{-1}]$,
	\item $R$ is an integral domain,
	\item $R/F$ is geometric,
	\item $R$ is $\partial$-simple.
\end{enumerate}
\end{defi}

We abbreviate the first condition by writing $R=F[W,W^{-1}]$. Clearly the field of fractions of a \PVR is
a \PVE. The next lemma implies that condition (2) is automatically satisfied, since $K$ will always be of
characteristic zero.

\begin{lem}\label{integral}
	Let $R$ be a $\partial$-simple $\partial$-ring and suppose that $R$ contains a field of
	characteristic zero. Then $R$ is an integral domain. 
\end{lem}
\begin{proof}
	Let $0 \ne f \in R$ be a zero-divisor. Then the localization map $R \to R_f$ is not
	injective and so $R_f = 0$ since $R$ is $\partial$-simple. This implies that $f$ vanishes on
	the prime spectrum of $R$ and must therefore be nilpotent.\\ 
	Consequently, the set of zero-divisors in $R$ coincides with the radical ideal $\sqrt{R}$
	of $R$. We claim that the latter ideal is a $\partial$-ideal. Indeed, suppose $f \in \sqrt{R}$ and
	let $n$ be minimal such that $f^n = 0$. Then we compute
	\[
		0 = \partial(f^n) = n f^{n-1} \partial(f) \text{.}
	\]
	But since $n$ was chosen minimal and $\Q \hra R$, we conclude that $\partial(f)$ is a
	zero-divisor. The first part of the proof yields $\partial(f) \in \sqrt{R}$. Clearly,
	$\sqrt{R}$ is proper since $1$ is not nilpotent and so $\sqrt{R} = 0$, because $R$ is
	$\partial$-simple. So $R$ is a domain, again by the first part of the proof.
\end{proof}

We cite the following classical existence result whose proof can be found in \cite{putsinger}.

\begin{theo}\label{existence} 
Let $F$ be a $\partial$-field with algebraically closed field of constants and let $A
\in \M_n(F)$. Then there exists a \PVR for $A$.
\end{theo}

We will deduce all fundamental results of Picard-Vessiot theory more or less directly from the
following key lemma. It generalizes Theorem 1 in \cite{levelt}. 

\begin{lem}\label{fundamental}
	Let $R = F[Y,Y^{-1}]$ be a \PVR over a $\partial$-field $F$ for a matrix $A \in \M_n(F)$. 
	Let $K$ denote the field of constants of $F$. Let $S/F$ be a geometric $\partial$-ring extension. 
	Assume there exists a fundamental solution matrix $X \in \GL_n(S)$ for $A$.\\
	Define $U$ to be the $K$-algebra of constants in $S \otimes_F R$. Then $U$ is a finitely
	generated $K$-algebra and the map
	\[
	\theta: S \otimes_K U \lra S \otimes_F R,\; s \otimes u \mapsto (s \otimes 1) u
	\]
	is an $S$-linear isomorphism of $\partial$-rings.
\end{lem}
\begin{proof} 
	In fact, the matrix $Z = X^{-1} \otimes Y \in \GL_n(S \otimes_F R)$ has coefficients in $U$
	since its derivative is the zero matrix.
	The surjectivity of $\theta$ is clear, since we have the equation
	\begin{equation}\label{surj}
	1 \otimes Y = (X \otimes 1) Z
	\end{equation}
	and $R$ is generated by $Y$.\\ 
	The injectivity is a little more difficult. Let $P$ be the kernel of the map $\theta$. 
	It is a $\partial$-ideal of 
	$S \otimes_K U$ and we will show that it is trivial. Choose a $K$-basis $\{ e_i \}$ of $U$
	where $i$ runs over some index set. Note that the elements $e_i$ are constants since $U$ consists of
	constants. Let $f = \sum s_i \otimes_K e_i \in P$ with $s_i \in S$ 
	such that the number of nonzero $s_i$ occuring in the sum is minimal. Call this number the
	length of $f$. If $s_i = 0$ for all
	$i$ we are done, so assume without loss of generality that $s = s_1 \ne 0$. The element $\partial(f)s -
	f\partial(s)$ lies in $P$ and has smaller length than $f$. So all its coefficients
	$\partial(s_i)s - s_i \partial(s)$ must vanish. Therefore the elements
	$\frac{s_i}{s}$ in the field of fractions of $S$ are constants so by assumption
	they lie in $K$. Thus we conclude that $f$ is of the form $s
	\otimes_K u$ for $u \in U$. Then, because $f$ lies in the kernel of
	$\theta$, we conclude $(s\otimes 1) u = 0$. But multiplication by $(s\otimes 1)$ in $S
	\otimes_F R$ is a flat base change of multiplication by $s \ne 0$ in $S$ and thus injective
	because $S$ is an integral domain. Thus $u = 0$ and finally $f = 0$.\\
	The restriction of $\theta$ to $K[Z,Z^{-1}] \subset U$ is surjective by (\ref{surj}) and 
	trivially injective. But then it is an isomorphism which implies $K[Z,Z^{-1}] = U$.
	Therefore, $U$ is a finitely generated $K$-algebra.
\end{proof}

\begin{cor}\label{pvr.inj.pve}
	Let $E/F$ be a Picard-Vessiot extension with algebraically closed field of constants $K$.
	Let $R/F$ be a \PVR for the same matrix. Then there
	exists an embedding 
	\[
	R \hra E 
	\]
	of differential rings over $F$.
	In particular, $E$ is the field of fractions of a Picard-Vessiot ring. If $E/F$ is generated
	by the fundamental matrix $W$ then its corresponding Picard-Vessiot ring is explicitly given
	by $F[W,W^{-1}]$.
\end{cor}
\begin{proof}
	In Lemma \ref{fundamental} set $S=E$. The restriction of the inverse map of $\theta$ to $1
	\otimes_F R$ gives an injection
	\[
	R \hra E \otimes_K U \text{.}
	\]
	Since $U$ is of finite type over $K$ and $K$ is algebraically closed, there exists a
	$K$-valued point $p \in \Hom_K(U, K)$. Composing the above injection with $1 \otimes p$
	yields the desired injection, because $R$ is $\partial$-simple.\\
	Let $E/F$ be a \PVE for a given matrix. 
	By \ref{existence}, there exists a Picard-Vessiot ring $R$ for the same matrix which can be
	embedded into $E$. But then clearly $E$ is the field of fractions of $R$.
	The last assertion is clear by comparing $W$ with the image of a generating fundamental 
	matrix of $R/F$ in $E$.
\end{proof}

One is tempted to prove the last assertion of the corollary directly, but the nontrivial assertion
here is the $\partial$-simplicity of $F[W,W^{-1}] \subset E$.\\
By using a little descent theory we obtain a result for non algebraically closed fields of
constants.

\begin{cor}\label{pve.pvr}
	Let $E/F$ be a \PVE generated by the fundamental matrix $W$. Then the ring $R = F[W,W^{-1}]
	\subset E$ is a \PVR over $F$. In particular, $E$ is the field of fractions of a \PVR. 
\end{cor}
\begin{proof}
	We have to show that $R$ is $\partial$-simple. By \cite[III.6.1]{stichtenoth} we know that 
	$E_{\Kb}/F_{\Kb}$ is a geometric $\partial$-field extension so it clearly is a \PVE. 
	By Corollary \ref{pvr.inj.pve} we conclude that $R_{\Kb} \subset E_{\Kb}$ 
	is $\partial$-simple. Suppose $P \subset R$ is a nonzero
	$\partial$-ideal. Because $\Kb/K$ is faithfully flat, $P_{\Kb}$ cannot be $0$ and thus
	$P_{\Kb} = R_{\Kb}$. So there exist $p_i \in P$ and $r_i \in R_{\Kb}$ such that $\sum r_i
	p_i = 1$. Summing over the finitely many Galois conjugates of both sides (with respect to
	$\Kb/K$) and using the assumption $\operatorname{char}(K) = 0$ implies $P = R$. 
	Thus $R$ is $\partial$-simple and therefore a \PVR over $F$.
\end{proof}

\begin{cor}\label{imply.dsimple}
	The first three conditions in Definition \ref{dgalois} imply the last condition.
\end{cor}
\begin{proof}
	Let $R/F$ be a $\partial$-ring extension satisfying the first three conditions of \ref{dgalois}. 
	Let $E/F$ be the corresponding extension of
	function fields. Clearly, $E/F$ is a \PVE. By Corollary \ref{pve.pvr} we obtain that $R/F$
	is a \PVR.
\end{proof}

The next theorem gives an existence result which is also valid in the case of a non algebraically 
closed field of constants.

\begin{theo}\label{existence2} 
	Let $S$ be a $\partial$-ring which is noetherian, integral and of finite Krull dimension. 
	Let $F$ be the field of fractions of $S$, whose field of constants we denote by $K$.
	Assume that $\Spec(S)$ has a $K$-valued point which is regular. 
	Then for any matrix $A \in \M_n(S) \subset \M_n(F)$ there exists a \PVR $R/F$.
\end{theo}
\begin{proof}
	Let $\O$ denote the structure sheaf of $Y = \Spec(S)$ and let $x \in Y(K)$ be the regular
	point which exists by assumption. 
	The stalk $\O_x$ is a regular local ring. Let $\m \subset \O_x$ denote the maximal
	ideal. The derivation on $S$ induces a derivation $\partial_x$ on $\O_x$ which can be
	identified with an element of the dual space of $\m/\m^2$. Therefore it is possible to 
	choose generators $t_1, \ldots, t_r$ of $\m$ reducing to a basis of $\m/\m^2$ such
	that $\partial_x(t_1) = 1 \text{(mod $\m^2$)}$ and $\partial_x(t_i) = 0 \text{(mod $\m^2$)}$ 
	for $i > 1$.\\
	Since $\O_x$ is regular we obtain
	\[
	\hat{\O}_x \cong K[ [t_1, t_2, \ldots, t_r]]
	\]
	where $\hat{\O}_x$ denotes the $\m$-adic completion of $\O_x$. We have injections
	\[
	S \hra \O_x \hra \hat{\O}_x\text{.}
	\]
	Recursively, we may replace $t_i$ by elements in $t_i + \m^2$ such that, letting $t=t_1$,  the
	derivation $\frac{d}{dt}$ on $K[ [t, t_2, \ldots, t_r]]$ induces the given 
	derivations on $\O_x$ and
	$S$. Clearly, the transformed parameters still generate $\hat{\O}_x$ as a power series
	ring. So we may assume $\partial = \frac{d}{dt}$.\\
	We want to construct a fundamental matrix $W \in \GL_n(\hat{\O}_x)$.
	Considering the fact that
	\[
	K[ [t, t_2, \ldots, t_r]] \cong K[ [t]]\,[ [ t_2, \ldots, t_r]]
	\]
	we may write
	\[
	A = A_0 + A_1 t + A_2 t^2 + \dots\quad\text{, $A_i \in K[ [ t_2, \ldots, t_r]]$}
	\]
	and make the Ansatz
	\[
	W = W_0 + W_1 t + W_2 t^2 + \dots\quad \text{, $W_i \in K[ [ t_2, \ldots, t_r]]$.}
	\]
	The desired condition $\partial(W) = AW$ yields the recursive equations
	\begin{align}\label{recursive}
		i W_{i} = \sum_{j+k = i-1} A_j W_k
	\end{align}
	for each $i >0$. As an initial condition we set $W_0$ to be the identity matrix and obtain a
	unique solution $W$ to (\ref{recursive}).\\
	We define $R = F[W,W^{-1}] \subset \Quot(\hat{\O}_x)$. Clearly, $R/F$ is geometric and
	integral, so it is a \PVR by \ref{imply.dsimple}.
\end{proof}

\begin{cor}
Let $K$ be a field and $A \in \M_n(\Kt)$. Then there exists a \PVE over $\Kt$ for $A$.
\end{cor}
\begin{proof}
Choose a point $x \in \P(K)$ such that $A \in \M_n(K[t]_x)$ (i.e. $x$ is a regular point of $A$). Now apply the 
preceding theorem to $S = K[t]_x$.
\end{proof}

The next proposition is a special case of Lemma \ref{fundamental}.
	
\begin{prop}\label{key}
	Let $R/F$, $R'/F$ be Picard-Vessiot rings for the same equation. Then the map 
	\[
	\theta: R' \otimes_K U \lra R' \otimes_F R, s \otimes u \mapsto (s \otimes 1) u
	\]
	is an $R'$-linear isomorphism of $\partial$-rings.
\end{prop}

Using this proposition, we deduce the representability of the functor $\underline{\Aut}^\partial(R/F)$
for a \PVR $R/F$. We actually prove a more general statement which also
implies the uniqueness of a \PVR for a given equation in the case of an algebraically closed field of 
constants.

\begin{cor}\label{isom.rep} 
	Let $R/F$, $R'/F$ be Picard-Vessiot rings for the same equation. Then the functor
	\[
	\underline{\Isom}^\partial(R,R') : ({\sf Algebras }/K) \lra ({\sf Sets}),\; L \mapsto
	\Isom^\partial(R_L, R'_L)
	\]
	is represented by $\Spec(U)$ over $K$. The isomorphisms are considered in $(\D/F_L)$ and $L$ is
	provided with the trivial derivation.
\end{cor}
\begin{proof}
	We have to show that $\Hom_K(U, -) \cong \underline{\Isom}^\partial(R,R')$. Let $L/K$ be a
	$K$-algebra. Then
	\begin{align*}
		\Hom_K(U,L) & \cong \Hom_{R'}^\partial(R' \otimes_K U, R'_L)\\
		& \cong \Hom_{R'}^\partial(R' \otimes_F R, R'_L)\\
		& \cong \Hom_F^\partial(R, R'_L)\\
		& \cong \Hom_{F_L}^\partial(R_L, R'_L)\\
		& \cong \Isom_{F_L}^\partial(R_L, R'_L)\text{,}
	\end{align*}
	where the rigidity in the last step can be concluded as follows:
	From the $\partial$-simplicity of $R$ one concludes that any element $f$ of
	$\Hom_F^\partial(R, R'_L)$ must be injective. But then a fundamental matrix generating $R/F$
	is mapped to a fundamental matrix $Y$ in $R'_L/F_L$. Observing that the constants in
	$R'_L$ equal $L$, this implies surjectivity of $f_L$ (Let $Z$ be a generating fundamental
	matrix of $R'_L/F_L$. Then $ZY^{-1}$ is a constant matrix and therefore $Y$ generates
	$R'_L/F_L$).
\end{proof}

Of course, this implies the representability of $\underline{\Aut}^\partial(R/F)$:
	
\begin{cor}\label{aut.rep}
	For a \PVR $R/F$ the group functor $\underline{\Aut}^\partial(R/F)$ 
	is represented by $G = \Spec(U)$, where $U$ is the $K$-subalgebra of constants in $R
	\otimes_F R$. Thus, $G \cong \underline{\Aut}^\partial(R/F)$ is a linear algebraic group over $K$.
\end{cor}

Here, by a linear algebraic group over $K$ we mean an affine group scheme of finite type over $K$.
We call $G$ the \defind{$\partial$-Galois group of $R/F$} and denote it by $\Gal(R/F)$.

\begin{cor}\label{storsor}
	Let $R/F$, $R'/F$ be Picard-Vessiot rings for the same equation and let $G =
	\Gal(R/F)$. Then the scheme $\underline{\Isom}^\partial(R,R')$ is a $G$-torsor over $K$.
	In particular, if $K$ is algebraically closed then a \PVR 
	for a given equation is unique up to isomorphism in $(\D/F)$.
\end{cor}
\begin{proof}
	The first assertion is immediate from Corollary \ref{isom.rep}. The second statement
	follows from the existence of a
	$K$-valued point of the scheme $\underline{\Isom}^\partial(R',R)$ which naturally
	translates into an isomorphism of $R'$ and $R$.
\end{proof}

To justify the introduction of Picard-Vessiot rings, we claimed they were geometric objects. 
This statement has in fact a very concrete meaning which we exhibit in the next theorem. It is a
generalization of the Torsor Theorem in \cite{putsinger}, though our proof will be different.
Indeed, all we have to do is to interpret the isomorphism in \ref{key} appropriately.

\begin{cor}[Torsor Theorem]\label{torsor}
	Let $R/F$ be a \PVR. Set $X = \Spec(R)$ and $Y = \Spec(F)$. 
	Then $X$ is a $G_Y$-torsor over $Y$.
\end{cor}
\begin{proof}
	The functor $\underline{\Aut}^\partial(R/F)$ acts on $X$. Under the isomorphism
	\begin{equation}\label{repisom}
		G \cong \underline{\Aut}^\partial(R/F)
	\end{equation}
	given explicitly in the proof of Corollary \ref{isom.rep} this action translates into
	an action of $G$ on $X$ (over $K$).\\
	Using the notation of \ref{key} (with $R' = R$), we claim that this action is given by the spectrum of
	\[
		\Delta: R \lra R \otimes_K U (= R \otimes_F U_F)
	\]
	where $\Delta$ is the restriction of the inverse of $\theta$ to $1 \otimes_F R$.\\
	Let $\phi \in \Aut^\partial(R_L/F_L)$ for a $K$-algebra $L$. Under the isomorphism (\ref{repisom}),
	$\phi$ is mapped to the restriction of
	\begin{equation}\label{explicit}
	\xymatrix{ R_L \otimes_{L} U_L \ar[r]^{\theta_L} & 
			R_L \otimes_{F_L} R_L \ar[r]^-{1 \otimes \phi} & R_L }
	\end{equation}
	to $1 \otimes_L U_L$, which yields an element $g_{\phi}$ in $\Hom_L(U_L, L)$. All we have to
	show is that we obtain $\phi$ if we compose $\Delta_L$ and $1 \otimes g_{\phi}$. But this 
	follows immediately from the explicit description of $g_{\phi}$ in (\ref{explicit}).\\
	Therefore $X$ is a $G_Y$-space with action given by
	\[
	\Spec(\Delta) : X \times_Y G_Y \lra X\text{.}
	\]
	Finally, Proposition \ref{key} implies that the map 
	\[
	X \times_Y G_Y \lra X \times_Y X
	\]
	induced by $\Spec(\Delta)$ and the first projection is an isomorphism or, 
	in other words, $X$ is a $G_Y$-torsor over $Y$.
\end{proof}

The next corollary generalizes parts of Corollary 1.30 in \cite{putsinger}.

\begin{cor} Let $R/F$ be a \PVR with field of fractions $E/F$. Let $X \to Y$
	be the spectrum of $F \hra R$ and let $G$ be the $\partial$-Galois group of $R/F$. Then
	\begin{enumerate}
		\item $X \to Y$ is smooth.
		\item The transcendence degree of $E$ over $F$ equals the dimension of $G$.
	\end{enumerate}
\end{cor}
\begin{proof}
	The smoothness follows from the Torsor Theorem by faithfully flat descent of this property.
	The transcendence degree of $E/F$ equals the relative dimension of $X \to Y$ which
	is equal to the dimension of $G$ by the Torsor Theorem.
\end{proof}

\section{Galois descent for Picard-Vessiot extensions}\label{descent}

Let $F$ be a $\partial$-field of characteristic zero and $\tF/F$ a $\partial$-field
extension which is an algebraic Galois extension of $F$.
Let $\Gamma$ be the Galois group of $\tF/F$. A \defind{descent datum} to $F$ in $(\D/\tF)$ is a pair
$(\tR, \rho)$ where $\tR$ is an object of $(\D/\tF)$ and $\rho:\Gamma \to \Aut^\partial(\tR)$ is a continuous action of
$\Gamma$ on $\tR$ such that the map $\tF \to \tR$ is $\Gamma$-equivariant.\\
The descent data form the objects of a category which we denote by $(\D/\tF)^\Gamma$ where the morphisms
are $\Gamma$-equivariant morphisms in $(\D/\tF)$.\\
Applying the base change functor to an object $R$ in $(\D/F)$ we obtain a descent datum. It is
easy to see that the resulting functor between $(\D/F)$ and $(\D/\tF)^\Gamma$ establishes an equivalence
of categories (this is an obvious extension of the classical result for vector spaces which in turn
is a trivial case of faithfully flat descent of quasi-coherent modules \cite[VIII]{sga1}). 
In other words, descent data to $F$ in $(\D/\tF)$ are effective.\\

We define the category $(\PV/F)$ of Picard-Vessiot extensions over $F$ as a full subcategory of
$(\D/F)$. Let $\tK/K$ be a Galois extension. The field of constants $K$ is algebraically closed 
in $F$, since any element which is algebraic over $K$ has to be constant itself. 
We conclude that $\tF := F \otimes_K \tK$ is a $\partial$-field with field of constants 
$\tK$ \cite[III.6.1]{stichtenoth}. Also, the extension
$\tF/F$ is Galois with Galois group isomorphic to the Galois group of $\tK/K$.\\ 
Let $(\PV/\tF)^\Gamma$ denote the full subcategory of $(\D/\tF)^\Gamma$ given by those descent data
whose underlying rings are Picard-Vessiot extensions of $\tF$.
The functor $-_{\tF}$ maps objects in $(\PV/F)$ to objects in $(\PV/{\tF})^\Gamma$.

\begin{theo}\label{pve.descent}
	The functor 
	\[
	{\mathcal B} : (\PV/F) \lra (\PV/\tF)^\Gamma
	\]
	which is given by base change is an equivalence of categories.
\end{theo}

\begin{proof}
	Let $(\tE, \rho)$ be an object in $(\PV/\tF)^\Gamma$. Descent data to $F$ in 
	the category $(\D/\tF)$ are effective and we obtain the corresponding $\partial$-field
	$E/F$ by taking invariants of $\tE$ under the action $\rho$. 
	The field of constants of $\tE$ being
	$\tK$ implies that the field of constants of $E$ is $K$ since the Galois action commutes with
	the derivation.
	Thus we only have to show that $E/F$ is generated by a fundamental matrix of a linear
	$\partial$-equation with entries in $F$.\\
	Let $\tE/\tF$ be generated by the fundamental matrix $W \in \GL_n(\tE)$ for the matrix
	$A \in \M_n(\tF)$. Since $\tE/E$ is Galois the field which is generated over
	$E$ by the entries of $W$ and $A$ is contained in a finite Galois extension $E'/E$. The
	field $E'$ is the field of invariants of $\tE$ under an open normal subgroup $\Delta \subset
	\Gamma$. We also define $F' := \tF^\Delta$ resp. $K' := \tK^\Delta$ obtaining finite Galois
	extensions of $F$ resp. $K$.\\
	Let $B = (1, b_2, \ldots, b_m)$ be a $K$-vector space basis of $K'$. It is also a basis of 
	$F'$ over $F$ and $E'$ over $E$, respectively.
	Note that the elements of $F' \subset E'$ are exactly the
	elements having coordinates in $F$ with respect to $B$.\\
	Consider the map
	\[
	\mu : E'  \lra \M_m(E)
	\]
	sending an element $x$ of $E'$ to the matrix representation of the endomorphism
	``multiplication by $x$'' with respect to the basis $B$. The image of $F'$ under
	$\mu$ lies in $\M_m(F)$. Since $B$ consists of constants the map $\mu$ commutes with the
	derivations and is thus a map of (noncommutative) $\partial$-rings (the derivation on 
	the ring $\M_m(E)$ being entry-wise).\\
	We compute
	\[
	\partial(\mu(W)) = \mu(\partial(W)) = \mu(AW) = \mu(A) \mu(W) \quad (\in \M_{mn}(E))
	\]
	where $\mu$ (and $\partial$) is applied to every entry of the respective matrices.
	Furthermore, we have
	\[
	1 = \mu(WW^{-1}) = \mu(W) \mu(W^{-1})
	\]
	and so $\mu(W) \in \GL_{mn}(E)$ is invertible 
	and thus a fundamental solution matrix of $\mu(A) \in \M_{mn}(F)$.\\
	All entries of $W$ are $K'$-linear combinations of entries of $\mu(W)$ thanks to the element
	$1$ being contained in the basis $B$. So $\tE/\tF$ is generated by the entries of $\mu(W)$.
	But this means that the inclusion map $F(\mu(W)_{ij}) \subset E$ becomes an isomorphism after
	applying the faithfully flat base change $- \otimes_F \tF$ and thus $F(\mu(W_{ij})) = E$.
\end{proof}

Let $E/F$ and $E'/F$ be Picard-Vessiot extensions and $K$ be the field of constants of $F$. For a field extension $L/K$
we call $E'$ an \defind{$L/K$-form} of $E$ if $E_L$ and $E'_L$ are isomorphic in $(\D/F_L)$.
	
\begin{cor}\label{forms.class}
	Let $E/F$ be a \PVE with $\partial$-Galois group $G$ and let $L/K$ be an algebraic 
	Galois extension with Galois group $\Gamma$. 
	Then the pointed Galois cohomology set $\Heins(\Gamma, G(L))$ classifies isomorphism classes of
	$L/K$-forms of $E/F$.
\end{cor}
\begin{proof}
This is clear since $G(L) \cong \Aut^\partial(E_L/F_L)$.
\end{proof}

The corollary shows that we have a correspondence between $L/K$-forms of a given \PVE $E/F$ and
$G$-torsors over $K$ which split over $L$ (where $G$ is the differential Galois group of $E/F$). 
On the level of Picard-Vessiot rings, we can make this correspondence explicit. 
Indeed, let $R \subset E$ be a \PVR (which exists by \ref{pve.pvr}).
To an $L/K$-form $E'/F$ of $E/F$ with corresponding \PVR $R' \subset E'$, we associate the functor 
$\underline{\Isom}^\partial(R,R')$ which is a $G$-torsor over $K$ by \ref{storsor}.\\
This result fits well with the theory of Tannakian categories where the fiber functors of a given
category are also classified by the torsors of the Tannakian fundamental group (cf.
\cite{deligne1}, \cite{deligne}). Actually, every $L/K$-form of a \PVE yields a fiber functor for a Tannakian category
associated to the original differential equation. A remarkable advantage over the general Tannakian
theory is, that we have an explicit description of the schemes which represent the 
$\underline{\Isom}$-functors in terms of the
differential-algebraic structure on the Picard-Vessiot rings.

\section{The Galois Correspondence}

The next theorem generalizes the Galois Correspondence to non algebraically closed fields of
constants. Using the theory developed in the last section it is obvious how to obtain a Galois 
correspondence by descending from the correspondence over $\Kb$.\\
However, we prefer to give a direct proof using the Torsor Theorem together with the theory of
descent for quasi-projective schemes.\\

First, we need a functorial version of invariants.
Let $S$ be a $K$-algebra. Let $H$ be a $K$-group functor which acts on $S$. By this we mean that for
any $K$-algebra $L$ the group $H(L)$ acts on $S_L$ and the action is functorial in $L$.
We say that $s \in S$ is \defind{invariant under $H$}, if for all $K$-algebras $L$ the element 
$s \otimes 1 \in S \otimes_K L$ is invariant under the action of the group $H(L)$. We write $S^H$ for the 
\defind{ring of invariants of $S$ under $H$}.

\begin{lem}\label{fields}
	Let $H$ be a linear algebraic group over $K$ acting functorially on a $K$-algebra $S$.
	Then $s \in S$ is invariant under $H$ if and only if it is invariant under $H(\Kb)$.
\end{lem}
\begin{proof}
	Let $H = \Spec(U)$. Then the action of the universal point in $H(U)$ given by $\id_U$
	determines a universal action $\Delta: S \to S \otimes_K U$.\\
	Suppose $s \ne 0$ is invariant under $H(\Kb)$. Extend $s_1 := s$ to a basis $\{s_i\}$ of $S$.
	Write $\Delta(s) = \sum_{i=1}^r s_i \otimes u_i$. By the invariance of $s$ we conclude that
	$u_i(g) = 0$ for $i > 1$ for every $g \in H(\Kb)$. So the functions $u_i$, $i>1$ vanish on the
	closed points of $H$ and by Hilbert's Nullstellensatz they vanish on all of $H$. Thus the
	$u_i$ are nilpotents but since any linear algebraic group in characteristic zero is reduced
	we conclude $u_i = 0$ for $i>1$. The same argument applied to $u_1 - 1$ shows that $u_1 = 1$
	,thus $\Delta(s) = s$ as we had to show.
\end{proof}

\begin{prop}\label{invariants}
	Let $F$ be a $\partial$-field with field of constants $K$. Let $E/F$ be a \PVE with
	$\partial$-Galois group $G$. Let $H$ be a closed $K$-subgroup of $G$. Then the inclusion $F
	\hra E^H$ is surjective if and only if $H = G$.
\end{prop}
\begin{proof}
	Let $R \subset E$ be a \PVR over $F$ (which exists by \ref{pve.pvr}). 
	By \ref{torsor}, $X = \Spec(R)$ over $Y = \Spec(F)$ is
	a $G$-torsor. Consider the functor which maps an $F$-algebra $S$ to the orbit set
	$X(S)/H_F(S)$ and denote its corresponding sheaf by $X/H$ (\cite[I, 5.5]{jantzen1}). 
	Since $X$ becomes isomorphic to the trivial torsor $G$ over a finite Galois extension of $F$, we
	conclude that $X/H$ becomes isomorphic to (the base change of) $G/H$. But by \cite[12.2.1]{springer}
	the functor $G/H$ is representable by a smooth quasi-projective scheme over $K$. So we may apply
	descent for quasi-projective schemes (\cite[VIII 7.7]{sga1}) to conclude that $X/H$ itself
	is representable by a smooth quasi-projective scheme over $F$. Since $X/H$ is a quotient of a
	connected scheme it must be connected and thus integral by smoothness. The map $X \to X/H$
	is an $H$-torsor (\cite[loc. cit.]{jantzen1}).\\
	The function field of $X/H$ is $E^H$ and the inclusion $F \hra E^H$ is the map on function fields 
	corresponding to the map $X/H \to Y$.\\
	In the case $H = G$, we obviously have $F = E^H$ since $X/G = Y$. Next assume 
	$F = E^H$. Then, on the one hand the map $X \to Y$ is the generic fiber of the map $X \to X/H$
	and thus an $H$-torsor, on the other hand $X \to Y$ is a $G$-torsor. But this implies that
	the inclusion $H \hra G$ becomes an isomorphism after the faithfully flat base change $X
	\times_K -$. So $H = G$.
\end{proof}

\begin{prop}\label{normal}
	Let $F$ be a $\partial$-field with field of constants $K$. Let $E/F$ be a \PVE with
	$\partial$-Galois group $G$. Let $H$ be a closed normal $K$-subgroup of $G$. 
	Then $E^H/F$ is a \PVE with $\partial$-Galois group $G/H$.
\end{prop}
\begin{proof}
As in the proof of Proposition \ref{invariants} we construct the scheme $X/H \to Y$ with with
function field $E^H$. Since $E^H \hra E$ it is clear that $E^H/F$ is geometric. By
\ref{imply.dsimple} it suffices to show that $E^H/F$ is generated by a fundamental matrix. 
By assumption $H \hra G$ is a normal subgroup. Then \cite[12.2.2]{springer} implies that $G/H$ is
indeed an affine $K$-group scheme. But then by descent for affine schemes \cite[VIII, 2]{sga1} the
scheme $X/H \to Y$ is affine itself. The natural action of $G/H$ turns $X/H$ into a $G/H$-torsor.\\
Choose a faithful linear representation $G/H \hra \GL_n$ over $K$. This defines a
left $G/H$-action on $\GL_n$. The fibration $X \times^G \GL_n$ is defined to be the quotient
of $X/H \times_K \GL_n$ under the right $G$-action given by $(x,y) \mapsto (xg,g^{-1}y)$. It
exists by \cite[I 5.14]{jantzen1} and is a $\GL_n$-torsor. Applying the flat base change $X/H \times_K -$ to the
injection $G/H \hra \GL_n$ we obtain the commutative diagram
\[
	\xymatrix@=4ex@C-1.2pc{
	X/H \times_K G/H\; \ar@{^{(}->}[rr] \ar[dr]& &X/H \times_K \GL_n \ar[dl]\\ 
	& X & 
	}
\]
which is $G/H$-equivariant. Therefore it descends to the quotients and yields
\begin{align}\label{tors.injection}
	\xymatrix@=4ex@C-1.2pc{
	X/H\; \ar@{^{(}->}[rr]\ar[dr] & & X/H \times^{G/H} \GL_n \ar[dl]\\ 
	& Y& 
	}
\end{align}
which is a $G/H$-equivariant embedding of $X/H$ into a $\GL_n$-torsor over $Y$.\\
By Hilbert 90 we have $\Heins(F, \GL_n) = 0$ and thus the latter $\GL_n$-torsor is trivial. So it is
isomorphic to the spectrum of the $F$-algebra $F[Z_{ij}, \det(Z)^{-1}] = F[Z,Z^{-1}]$. Suppose that
$X/H$ is the spectrum of the $F$-algebra $S$. Then the embedding (\ref{tors.injection}) induces a
surjection $F[Z,Z^{-1}] \to S$. So $S = F[W,W^{-1}]$ where $W$ denotes the image of the matrix $Z$ in $S$.
Since the map is $G/H$-equivariant, the $G/H$-action on $S$ is given by mapping $W$ to $WC$ for $C
\in (G/H)(L) \hra \GL_n(L)$ (for a $K$-algebra $L$). Thus the matrix $A := \partial(W)W^{-1}$ is
$G/H$-invariant and has therefore entries in $F = S^{G/H}$. So $W$ is a fundamental matrix for $A
\in \M_n(F)$ and generates $E^H$.\\
We have a natural inclusion $G/H \hra \Gal(E^H/F)$ which becomes an isomorphism after the
faithfully flat base change $X/H \times_K -$, so this inclusion must itself be an isomorphism.
\end{proof}

Now we are ready to prove the Galois Correspondence.

\begin{theo}\label{correspondence}
	Let $F$ be a differential field with field of
	constants $K$ and let $E$ be a \PVF over $F$ with $\partial$-Galois group $G$.
	Let
	\[
	{\mathcal H} = \{ H \le G | \; \text{$H$ closed $K$-subgroup scheme of $G$} \}
	\]
	as well as
	\[
	{\mathcal M} = \{ M | \; \text{$F \hra M \hra E$ extensions of differential fields } \}\text{.}
	\]
	\begin{enumerate}
		\item The map
			\[
			\Psi: {\mathcal M} \lra {\mathcal H},\; M \mapsto \Gal(E/M)
			\]
			is an anti-isomorphism of lattices.
        		The inverse map of $\Psi$ is given by
        		\[
        		\Phi: {\mathcal H} \lra {\mathcal M},\; H \mapsto E^H \text{.}
        		\]
		\item If $H \in {\mathfrak H}$ is a normal subgroup scheme of $G$ 
			then $E^H$ is a \PVE over $F$ with $\partial$-Galois group $G/H$.
	\end{enumerate}
\end{theo}

\begin{proof}
We proved (2) in \ref{normal}. For (1) we first have to show that $\Psi$ is well-defined. Let $M \in
\mathcal M$. We have a natural inclusion $\Gal(E/M) \hra \Gal(E/F)$ and we have to show that it is
closed. Let $R \hra E$ be a \PVR, then we showed in \ref{aut.rep} that $\Gal(E/F)$ is representable
by the spectrum of the $K$-algebra $(R\otimes_F R)^\partial$ where $-^\partial$ means that we
consider constants with respect to $\partial$. Thus, the group $\Gal(E/M)$ is represented by the
spectrum of $(R \otimes_M R)^\partial$. The inclusion $\Gal(E/M) \hra \Gal(E/F)$ is given by the
spectrum of the natural surjection $(R \otimes_M R)^\partial \to (R\otimes_F R)^\partial$ of $K$-algebras. 
The kernel of this surjection is the defining ideal of $\Gal(E/M)$ as a closed $K$-subgroup of
$\Gal(E/F)$.\\
From \ref{invariants}, we directly conclude that $\Psi$ and $\Phi$ are inverse to each other.
\end{proof}

\section{Examples} \label{sec.ex}
In this section we study some explicit examples.
We consider the case of the Galois extension $\C/\R$ with Galois group $\Gamma$
generated by the complex conjugation $\tau$. 

	\subsection{$G \cong \Gm$} 
		The differential equation $\frac{d}{dt}(y) = y$ over $\Ct$ has the ring $\tR =
		\Ct[y,y^{-1}]$ as \PVR. 
		We extend the complex conjugation to $\tR$ via $\tau(y) =
		y$. The ring of invariants under $\tau$ is given by $R = \Rt[y,y^{-1}]$ which is a
		\PVR by Theorem \ref{pve.descent}. The Galois group scheme $G$ of $R$ over $\Rt$ is easily
		seen to be the group scheme $\Gm$ over $\R$.\\
		By Hilbert 90, we know that $\Heins(\Gamma,\Gm(\C))$ is trivial and 
		conclude by Theorem \ref{forms.class} that $R$ is the only $\R$-structure on
		$\tR$. 
	\subsection{$G \cong \mu_3$} \label{example.mu}
		Let $\tilde{E} = \Quot(\tilde{R})$ with $\R$-structure $E = \Quot(R)$ from the
		previous example.
		Let $\mu_3$ be the closed normal subgroup scheme of $\Gm$ defined by the equation $z^3 =
		1$. We have $E^{\mu_3} = \Rt(y^3)$ corresponding to $\mu_3$ by the Galois
		correspondence. Note that we need to consider functorial invariants since $\mu_3(\R)
		= 1$. Also note that $E/\Rt(y^3)$ is a finite $\partial$-Galois extension but not
		a finite Galois extension in the usual sense. This comes from the fact that in the 
		theory described here, we also consider finite \'etale torsors as $\partial$-Galois
		(since the Galois Correspondence includes them).
	\subsection{$G \cong \SO_2$} \label{example.so}
		Consider the differential equation $\partial(y) = iy$ over $\Ct$ which induces the
		\PVR $\tilde{R} = \Ct[y,y^{-1}]$. There is an analytical interpretation of $y$ as
		the function $\exp(it)$. 
		We extend the complex conjugation to a
		differential automorphism of $\tR$ by letting $\tau(y) = y^{-1}$. The corresponding
		ring of $\Gamma$-invariants $R$ is a \PVR over $\Rt$ for the equation
		\[
		\partial(Y) = \mu(i) Y = \left( \begin{array}{cc} 0 & -1\\ 1 & 0 \end{array} \right) Y\text{.}
		\]
		Here we use the notation of the proof of Theorem \ref{pve.descent} where we choose
		the $\R$-basis $B = (1,i)$ of $\C$. The base extension $R \otimes -$ yields the $R$-basis
		$(1 \otimes 1, 1 \otimes i)$ of $\tR$.
		With respect to this basis, we write
		\[
		y = \underbrace{ \frac{1}{2}\left( y + \tau(y) \right) }_{\cos(t)} \otimes 1
		+ \underbrace{ \frac{1}{2i}\left( y - \tau(y) \right)}_{\sin(t)} \otimes i\text{.}
		\]
		Still following the proof of Theorem \ref{pve.descent}, the matrix
		\[
		\mu(y) = \left( \begin{array}{cc} \cos(t) & -\sin(t)\\ \sin(t) & \cos(t) \end{array} \right)
		\]
		is a generating fundamental matrix of $R = \Rt[\cos(t), \sin(t)]$ with the only relation
		$\cos(t)^2 + \sin(t)^2 = 1$. 
		The Galois group scheme of $R$ over $\Rt$ is the group scheme $\SO_2$ over $\R$. 
		By \cite[p.141]{serre}, there is a bijection between $\Heins(\Gamma,\SO_2(\C))$ and the set
		of classes of quadratic forms of rank $2$ which have positive discriminant. We
		conclude that $\Heins(\Gamma, \SO_2(\C))$ consists of two elements. A nontrivial
		cocycle $\chi$ is given by 
		\[
		\chi(\tau) = \left( \begin{array}{rr} -1 & 0\\ 0 & -1 \end{array} \right) \text{.}
		\]
		Twisting the above $\Gamma$-action ($\tau(y) = y^{-1}$) with $\chi$ we obtain the
		action $\tau_\chi(y) = -y^{-1}$. This action corresponds to the \PVR $R' = \Rt[i \cos(t), i \sin(t)]$ with
		the relation $(i\cos(t))^2 + (i\sin(t))^2 = -1$ and according to Theorem
		\ref{pve.descent}, the ring $R'$ is a nontrivial $\C/\R$-form of $R$. A fundamental
		matrix for $R'$ is given by
		\[
		\mu_\chi(y) = \left( \begin{array}{cc} i \sin(t) & -i \cos(t) \\ i \cos(t) & i \sin(t)
		\end{array} \right) \text{.}
		\]
		The Galois group scheme
		of $R'$ over $\Rt$ is again $\SO_2$. In general, different forms of a
		given \PVR do not need to have isomorphic Galois group schemes. 

	\subsection{$G \cong \C^*$ (as $\R$-structure on $\Gm^2$)} \label{example.C}
		We study the equation defined by
		\[
		A = \left( \begin{array}{rr} 1 & -1\\ 1 & 1 \end{array} \right)
		\]
		over $\Ct$. 
		A \PVR for $A$ over $\Ct$ is given by 
		\[
		\tR = \Ct[x,y,z]/( (x^2 + y^2)z - 1)
		\]
		generated by the fundamental matrix
		\[
		Y = \left( \begin{array}{rr} x & -y\\ y & x \end{array} \right) \text{.}
		\]
		An analytical interpretation would be
		\begin{align*}
			x = e^t \cos(t)\\
			y = e^t \sin(t)\\
			z = e^{-2t} \text{.}
		\end{align*}
		Taking invariants under the $\Gamma$-action given by $\tau(x) = x$ and $\tau(y) = y$
		yields the $\R$-structure 
		\[
		R = \Rt[x,y,z]/( (x^2+y^2)z - 1 )
		\]
		which has the Galois group scheme defined by the Hopf algebra
		\[
		U = \R[a,b,c]/ ( (a^2+b^2)c - 1 ) \text{.}
		\]
		Here the Hopf algebra structure is given by matrix operations on the universal matrix
		namely
		\begin{align*}
		\left( \begin{array}{rr} a & -b\\ b & a \end{array} \right) 
			\mapsto &
		\left( \begin{array}{rr} a & -b\\ b & a \end{array} \right) 
			\otimes
		\left( \begin{array}{rr} a & -b\\ b & a \end{array} \right) 
			 \\ = &
		\left( \begin{array}{rr} a \otimes a - b \otimes b & -a \otimes b - b \otimes a\\ 
				b \otimes a + a \otimes b & -b \otimes b + a \otimes a \end{array} 
				\right) 
		\end{align*}
		for comultiplication and similarly for coinverse and counit.
		Thus we see that $G$ is a $\C/\R$-form of $\Gm^2$ with
		\[
		G(\R) \cong \C^* \text{.}
		\]
		By \cite[12.3.4. Proposition]{springer}, there is an exact sequence of groups 
		(all group schemes are commutative)
		\begin{align*}
			1 & \lra \SO_2(\R) \lra G(\R) \overset{\det}{\lra} \Gm(\R) \lra\\
		&\overset{\delta}{\lra} \Heins(\Gamma, \SO_2(\C)) 
		\overset{1}{\lra} \Heins(\Gamma, G(\C)) \overset{j}{\lra} \Heins(\Gamma,\Gm(\C)) 
		\end{align*}
		where the ``1'' comes from the surjectivity of $\delta$. This implies that the map $j$
		is injective and thus $\Heins(\Gamma, G(\C))$ is trivial, since
		$\Heins(\Gamma,\Gm(\C))$ is trivial.
		Therefore, no nontrivial $\C/\R$-forms of $R$ exist.

\section{Monodromy and Galois actions}\label{monodromy}

We consider the noncompact Riemann surface $X$ which is obtained by removing a finite, nonemtpy set
$S$ of complex conjugate points from the Riemann sphere $\P$. We assume that $\infty$ is contained
in $S$. To emphasize the fact that we are working over the complex numbers, we switch from the
variable $t$ to the variable $z$.\\ 
Let $\Gamma$ denote the Galois group of $\C/\R$ with generator $\tau$. 
The group $\Gamma$ acts continuously on $X$ in the obvious way. We identify the sheaf of
holomorphic functions $\O$ on $X$ with its \'espace \'etal\'e $p: \O \to X$. So $\O$ is a Riemann
surface with underlying space $\coprod_{x\in X} \O_x$ and $p$ is a local homeomorphism (see
\cite[Theorem 6.8]{forster}).

\begin{prop}\label{cont} If $f$ is a holomorphic function on $U \subset X$ then the function 
	$\tau.f$ given by
	\[
	\tau \circ f \circ \tau^{-1}
	\]
	is homolorphic on $\tau(U)$. We have the formula
	\begin{equation}\label{cont.der}
		\pd{\tau.f}{z} = \tau.\pd{f}{z}\text{.}
	\end{equation}
	In particular, this induces a continuous action of the group $\Gamma$ on $\O$.
\end{prop}

\begin{proof}
	Chain rule.
\end{proof}

Next, we want to analyze how this action interacts with analytic continuation. Assume that $f_a$ is
a germ of a holomorphic function at the point $a \in X$. By the proposition, $\tau.f_a$ is a germ of
a holomorphic function at $\bar{a}$. Let $I$ be the real unit interval and 
\[
\alpha : I \lra X
\]
be a closed path with $\alpha(0) = a$. Let $\tilde{\alpha}$ be a lifting of $\alpha$ to $\O$ with
$\tilde{\alpha}(0) = f_a$, i.e. an analytic continuation of $f_a$ along $\alpha$. Then we define
\[
(\tau.\tilde{\alpha})(t) := \tau.(\tilde{\alpha}(t))\text{.}
\]
Because the action of $\Gamma$ is continuous this yields a closed path in $\O$ with
\[
(\tau.\tilde{\alpha})(0) = \tau.f_a
\]
which lifts the path $\tau \circ \alpha$.
Thus we conclude that $\tau.\tilde{\alpha}$ is the (unique) analytic continuation
of $\tau.f_a$ along $\tau \circ \alpha$.\\

This observation has implications on the monodromy representation of a linear differential
equation. Namely let $A \in \M_n(\Cz)$ with complex conjugate singularities $S$ 
and let $a \in X(\R)$ be a real regular point, meaning that all the entries of $A$ are holomorphic
at $a$. The general theory ensures
the existence of a fundamental solution matrix $Y \in \GL_n(\O_a)$, i.e.
\[
\pd{}{z}Y = AY\text{.}
\]
The monodromy representation
\[
\mu_Y : \pi_1(X,a) \lra \GL_n(\C)
\]
is given by analytic continuation of $Y$ along paths in $X$ with base point $a$. The image of this
homomorphism is called the \defind{monodromy group} of $A$. If we replace $Y$
by another fundamental matrix for $A$ then we obtain an equivalent representation. Thus the
monodromy group is defined up to conjugation.\\

\begin{prop}\label{mon.conj}
	If $\mu_Y$ is the monodromy representation of $A \in \M_n(\Cz)$ then the monodromy
	representation of $\tau.A$ induced by $\tau.Y$ is given by
	\[
	\pi_1(X,a) \lra \GL_n(\C), \alpha \mapsto
	\tau.\mu_Y(\tau^{-1} \circ \alpha)\text{.}
	\] 
\end{prop}
\begin{proof}
	Using formula (\ref{cont.der}) from Proposition \ref{cont} we obtain
	\[
	\pd{}{z}(\tau.Y) = \tau.(\pd{}{z}Y) = \tau.(AY) = \tau.A \tau.Y
	\]
	thus $\tau.Y$ is a fundamental solution matrix of $\tau.A$ at $\tau(a) = a$.
	Let $\alpha$ be a loop with base point $a$. By above remarks, the analytic continuation of
	$\tau.Y$ along $\tau \circ \alpha$ is given by applying $\tau$ to the analytic continuation
	of $Y$ along $\alpha$. Therefore 
	\begin{align*}
		\tau.Y \mu_{\tau.Y}(\tau \circ \alpha) &= \tau.(Y \mu_Y(\alpha))\\
		&= \tau.Y \tau.\mu_Y(\alpha)\text{,}
	\end{align*}
	which implies
	\[
	\mu_{\tau.Y}(\tau \circ \alpha) = \tau.\mu_Y(\alpha) \text{.}
	\]
	The result follows by replacing $\alpha$ with $\tau^{-1} \circ \alpha$.
\end{proof}

\section{Descent with Riemann-Hilbert}\label{riemann}

The following corollary together with the Riemann-Hilbert correspondence will be the key to solving
the inverse problem over $\Rz$. We say that two $\partial$-equations $\pd{y}{z} = Ay$ and $\pd{y}{z} = By$
with $A, B \in \M_n(\Cz)$ are \defind{equivalent} if there exists $C \in \GL_n(\Cz)$ such that
$C^{-1} A C - C^{-1} \pd{C}{z} = B$.

\begin{prop}\label{dipl}
	Let $\tilde{E}/\Cz$ be a \PVE for the matrix $A \in \M_n(\Cz)$. 
	If the equation $\pd{y}{z} = Ay$ is equivalent to its conjugate
	$\pd{y}{z} = (\tau.A) y$, then $\tilde{E}/\Cz$ descends to a \PVE $E/\Rz$.
\end{prop}
\begin{proof}
	Without restriction let $0$ be a regular point for $A$. Then we obtain a fundamental solution
	matrix $Y \in \GL_n(\C(( z )))$ with entries in the field of Laurent series. Here the
	derivation is extended to the field $\C(( z ))$ in the obvious way such that the field of
	constants is still $\C$. By definition the field $\Cz(Y) \subset \C(( z ))$ is a \PVE
	over $\Cz$. We may denote it by $\tE$ because Picard-Vessiot extensions are unique over $\Cz$.
	The group $\Gamma=\operatorname{Gal}(\C/\R)$ acts on $\C(( z ))$ by acting on the 
	coefficients of the power series.\\
	The asumption that $A$ is equivalent to $\tau(A)$ implies that this 
	Galois action restricts to an action on $\tE$. Thus \ref{pve.descent} implies the assertion.\\
\end{proof}

We remark that one can generalize this result to descend from $\overline{K}(z)$ to $K(z)$ for any field $K$.\\
Let $\operatorname{RegSing}(\Cz, S)$
be the category of regular singular linear differential equations over $\Cz$ with singularities contained in the
finite set $S \subset \P$. Let $\operatorname{\sf Repr}_{\pi_1}$ be the category of finite-dimensional 
complex linear representations of the fundamental group $\pi_1(\P \backslash S)$. 
Both categories are neutral Tannakian categories over $\C$.\\
The main ingredient of the proof will be the Riemann-Hilbert correspondence which is stated in the
following theorem. 
A proof may be found in \cite[Theorem 6.15]{putsinger}.

\begin{theo}\label{mon}
	The functor 
	\[
	\mathcal{M} : \operatorname{\sf RegSing}(\Cz, S) \lra \operatorname{\sf Repr}_{\pi_1}
	\]
	which associates to an equation its monodromy representation is an equivalence of Tannakian
	categories.
\end{theo}

This theorem contains two striking facts both of which we will use in the proof of the main theorem.
The first one is the existence of a regular singular equation for any prescribed monodromy
representation.
Secondly, a regular singular differential equation is determined up to equivalence by its
monodromy representation: 

\begin{cor}\label{mon.equ}
	Two regular singular equations over $\Cz$ are equivalent if and only if their monodromy
	representations are equivalent.
\end{cor}

\section{The inverse problem over $\Rz$}\label{finalsection}

Now we combine the above results to prove the main theorem. Using descent theory one
obtains an equivalence of categories between linear algebraic groups over $\R$ and such groups over
$\C$ together with a (compatible) action of the Galois group $\Gamma = \operatorname{Gal}(\C/\R)$
(\cite[17.3]{waterhouse}).

\begin{theo}\label{inverse}
	Every linear algebraic group $G$ over $\R$ is the differential Galois group of a regular
	singular differential equation over $\Rz$ .
\end{theo}
\begin{proof}
	Suppose $G \subset \GL_{n,\R}$ as group schemes over $\R$. 
	Then $G(\R)$ is the group of elements of $G(\C) \subset \GL_n(\C)$ which are invariant under 
	the natural action of $\Gamma$. By Lemma 5.13 in \cite{putsinger} there exist matrices $C_1,
	\ldots, C_r$ which generate a Zariski-dense subgroup of $G(\C)$.\\
	We choose $r$ points in the complex upper halfplane and set
	\[
	S := \left\{ s_1, \ldots, s_r, \tau(s_1), \ldots, \tau(s_r), \infty \right\}\text{.}
	\]
	Let $X$ denote the complement of $S$ in $\P$. We choose a real basepoint $a \in \R$.
	For $1 \le i \le r$ let $\alpha_i$ be a loop around $s_i$ which does not enclose any $s_j$ 
	for $j \ne i$. Then the group $\pi_1(X,a)$ is freely generated by
	the homotopy classes of $\alpha_1, \ldots, \alpha_r, \tau\circ \alpha_1, \ldots,
	\tau\circ\alpha_r$.\\
	By Theorem \ref{mon} there exists a regular singular equation given by $A \in \M_n(\Cz)$
	whose monodromy representation has the properties
	\[
		\mu_Y(\alpha_i) = C_i  \quad\text{and}\quad
		\mu_Y(\tau \circ \alpha_i) = \tau(C_i)\text{.}
	\]
	The differential Galois group of $A$ over $\Cz$ is the group $G(\C)$ since for regular
	singular differential equations the monodromy group is Zariski-dense in the 
	differential Galois group (\cite[Theorem 5.8]{putsinger}).\\
	Using Proposition \ref{mon.conj}, we compute the monodromy represention of the equation defined by
	$\tau.A$: 
	\begin{align*}
		\mu_{\tau.Y}(\alpha_i) &= \tau.(\mu_Y(\tau^{-1} \circ \alpha_i))\\
				&= C_i
	\end{align*}
	and analogously
	\[
		\mu_{\tau.Y}(\tau \circ \alpha_i) = \tau(C_i)\text{.}
	\]
	Thus the monodromy representations of $A$ and $\tau.A$ are equivalent and Corollary
	\ref{mon.equ} implies that the equations defined by $A$ and $\tau.A$ correspondingly must be
	equivalent.\\
	Let $\tilde{E} \subset \O_a$ be the \PVE generated over $\Cz$ by the above matrix $Y \in
	\GL_n(\O_a)$. Then Proposition \ref{dipl} and its proof imply the existence of an action of
	$\Gamma$ on $\tilde{E}$ commuting with the derivation and extending the action on $\Cz$. 
	The real structure $E/\Rz$ of $\tilde{E}/\Cz$ is obtained by taking invariants under this
	action.\\
	The induced action on the differential Galois group $G(\C) = \Aut^\partial(\tilde{E}/\Cz)$ is then
	given by mapping $\phi \in G(\C)$ to $\tau \circ \phi \circ \tau^{-1}$.
	Now let $\phi$ be the automorphism induced by analytic continuation along a path $\alpha$,
	i.e. $\phi(Y) = Y \mu_Y(\alpha)$ as well as $\phi(\tau.Y) = \tau.Y \mu_{\tau.Y}(\alpha)$. Then
	\begin{align*}
	\tau \circ \phi \circ \tau^{-1} (Y) &= \tau.(\phi(\tau.Y))\\
					&= \tau.(\tau.Y\mu_{\tau.Y}(\alpha))\\
					&= Y \tau.(\mu_Y(\alpha))\text{,}
	\end{align*}
	where the last step follows, because $\mu_{\tau.Y}(\alpha) = \mu_Y(\alpha)$ holds due to
	above calculations.\\
	So on the monodromy group, 
	the action of $\Gamma$ on $G(\C) \subset \GL_n(\C)$ coincides 
	with the action defining the real algebraic group $G$. But since the monodromy group is
	Zariski-dense inside $G(\C)$, the action must coincide on the whole group $G(\C)$.\\
	Therefore the differential Galois group of $E/\Rz$ is $G$.
\end{proof}

\bibliographystyle{amsalpha}
\protect
\bibliography{real_inverse}

\end{document}